\documentclass[11pt]{amsart} 
\usepackage{fancyhdr}
\usepackage{enumitem}
\usepackage{amsmath,amsthm,amssymb,mathrsfs,url}
\usepackage{tikz-cd}
\usepackage{marginnote}
\newcommand{\filename}{Abstract-Simplicial-Complexes-13-March-2025.tex} 

\providecommand{\binom}[2]{{#1\choose#2}}

\renewcommand{\geq}{\geqslant}
\renewcommand{\leq}{\leqslant}

\newcommand{\PP}{\mathbb{P}} 
\newcommand{\RR}{\mathbb{R}} 
\newcommand{\ZZ}{\mathbb{Z}} 

\newtheorem{theorem}{Theorem}[section]

\theoremstyle{definition}

\newtheorem{example}[theorem]{Example}

\numberwithin{equation}{section}

\title{Abstract simplicial complexes in {\tt Macaulay2}}

\author{Nathan Grieve}

\address{Department of Mathematics \& Statistics,
Acadia University, Huggins Science Hall, Room 130,
12 University Avenue,
Wolfville, NS, B4P 2R6
Canada; School of Mathematics and Statistics, Carleton University, 4302 Herzberg Laboratories, 1125 Colonel By Drive, Ottawa, ON, K1S 5B6, Canada; 
D\'{e}partement de math\'{e}matiques, Universit\'{e} du Qu\'{e}bec \`a Montr\'{e}al, Local PK-5151, 201 Avenue du Pr\'{e}sident-Kennedy, Montr\'{e}al, QC, H2X 3Y7, Canada; Department of Pure Mathematics, University of Waterloo, 200 University Avenue West, Waterloo, ON, N2L 3G1, Canada
}
\email{nathan.m.grieve@gmail.com}%

\begin{document}

\begin{abstract}
{\tt AbstractSimplicialComplexes.m2} is a computer algebra package written for the computer algebra system {\tt Macaulay2} \cite{M2}.  It provides new infrastructure to work with abstract simplicial complexes and related homological constructions.  Its key novel feature is to implement each given abstract simplicial complex as a certain graded list in the form of a hash table with integer keys.    Among other features, this allows for a direct implementation of the associated reduced and non-reduced simplicial chain complexes.  Further, it facilitates construction of random simplicial complexes.  The approach that we employ here builds on the {\tt Macaulay2} package {\tt Complexes.m2} \cite{Stillman:Smith:Complexes.m2}.  It complements and is entirely different from the existing {\tt Macaulay2} simplicial complexes framework that is made possible by the package {\tt SimplicialComplexes.m2} \cite{Smith:et:al:SimplicialComplexes.m2:jsag}.  
\end{abstract}

\thanks{
\emph{Mathematics Subject Classification (2020):}  05E45, 55U10, 62R40, 19-04. \\
\emph{Key Words:} Simplicial complexes, Simplicial chain complexes, Random simplicial complexes, Topological Data Analysis. \\
The author thanks the Natural Sciences and Engineering Research Council of Canada for their support through his grants DGECR-2021-00218 and RGPIN-2021-03821. \\
Date: \today.  \\
File name: \filename}

\maketitle

\section{Introduction}\label{intro}

{\tt AbstractSimplicialComplexes.m2} is a software package for {\tt Macaulay2 (version 1.24.11)} \cite{M2}.  It  provides new infrastructure to work with abstract simplicial complexes.  Its key novel feature is to implement each given abstract simplicial complex as a graded list (in the form of a hash table with integer keys).  This approach is entirely different from the methodology used in the existing {\tt Macaulay2} package {\tt SimplicialComplexes.m2} \cite{Smith:et:al:SimplicialComplexes.m2}, \cite{Smith:et:al:SimplicialComplexes.m2:jsag}. 

The approach that we employ here builds on Stillman and Smith's package {\tt Complexes.m2} \cite{Stillman:Smith:Complexes.m2}.  Among other features, our approach here allows for a direct implementation of simplicial and reduced simplicial chain complexes.  Another method which is of an independent interest  produces as output the  induced chain complex maps that arise via subsimplicial complexes.  Further we provide methods for producing random simplicial complexes and, in particular, this allows for calculation of random simplicial homology complexes. (See Section \ref{random:simplicial:complexes}.)

To place matters into perspective, while the existing implementation of simplicial complexes within {\tt SimplicialComplexes.m2} is well suited for simplicial complex commutative algebra questions that arise via the celebrated \emph{Stanley-Reisner correspondence}, a key limitation occurs when working with simplicial complexes and related homological questions.  Indeed, when working within {\tt SimplicialComplexes.m2}  there is the need to introduce a suitable polynomial ring and then to regard the vertices of each given simplicial complex as variables of that polynomial ring.
  
By contrast, {\tt AbstractSimplicialComplexes.m2} implements a framework for simplicial complexes and related homological questions that is based on graded lists as a primitive data type and without need to introduce auxiliary constructs such as polynomial rings.  (We refer to \cite{CCA} or \cite{Bruns-Herzog}, for example, for the most basic details about the Stanley-Reisner correspondence and to texts such as \cite{Armstrong:basic:topology}, \cite{Hatcher:Alg:Top} and \cite{Munkres:Alg:Top} for the  most basic topological theory of abstract  simplicial complexes.)

The advantages and drawbacks to such design choice considerations---depending on the particular applications that a given user has in mind---have indeed been well known to {\tt Macaulay2} developers for some time.  For example, they were highlighted at early stages in the development of the {\tt Macaulay2} package {\tt SpectralSequences.m2} \cite{M2:Spectral:Sequences}. 

Indeed, the primary impetus for what we do here arose after revisiting the design choice decisions that were made in creating the package {\tt SpectralSequences.m2}.  The aim was to extend its functionality and to give a good forward compatible integration with the package {\tt Complexes.m2}.   An overview of the main features of the package {\tt SpectralSequences.m2} is provided in  \cite{Boocher:Grieve:Grifo}.  Such decision choice considerations for the package {\tt AbstractSimplicialComplexes.m2}, which are compatible with the package {\tt Complexes.m2}, are a key novel feature to our present work here.

There are also pedagogical reasons for being able to interact with computational tools for understanding homological questions related to simplicial complexes that are independent of the correspondence of Stanley and Reisner.  This is illustrated in Example \ref{Eg:PP^2:RR:Klein:bottle} below.

\begin{example}[Calculating the homology groups of $\PP^2_{\RR}$ and the Klein bottle]\label{Eg:PP^2:RR:Klein:bottle}
The following code snippets illustrate how to calculate the non-reduced simplicial homology groups of $\PP^2_{\RR}$ and the Klein bottle.  Related calculations are illustrated within the packages  {\tt SpectralSequences.m2} \cite{M2:Spectral:Sequences} and {\tt SimplicialComplexes.m2} \cite{Smith:et:al:SimplicialComplexes.m2}.  The simplicial complex realizations that we use here are compatible with what is done in \cite{M2:Spectral:Sequences}; they follow the approach of the text \cite{Armstrong:basic:topology}.  The calculations that we do here should be compared with the approach to doing them within the framework of {\tt SimplicialComplexes.m2} \cite{Smith:et:al:SimplicialComplexes.m2}.  The simplicial complex conventions used in \cite{Smith:et:al:SimplicialComplexes.m2} are based on those of \cite{Munkres:Alg:Top}.  We refrain from producing such comparable calculations here; instead we focus on illustrating how to perform such calculations using the package {\tt AbstractSimplicialComplexes.m2}.
\begin{scriptsize}
\begin{verbatim}
nathangrieve@Nathans-Air-3 ~ % M2
Macaulay2, version 1.24.11
with packages: ConwayPolynomials, Elimination, IntegralClosure, InverseSystems, 
Isomorphism, LLLBases, MinimalPrimes, OnlineLookup,
PrimaryDecomposition, ReesAlgebra, Saturation, TangentCone, Truncations, Varieties
i1 : needsPackage"AbstractSimplicialComplexes";
i2 : -- A simplicial complex realization of PP^2_RR 
as given in Armstrong's book "Basic Topology"
     K = abstractSimplicialComplex {{1,2,3},{2,3,4},{3,4,5}, {1,5,4},
     {5,2,1},{5,6,2},{4,6,2},
     {1,6,4},{3,6,5},{1,6,3}};
i3 : prune HH simplicialChainComplex K
       1
o3 = ZZ  <-- cokernel | 2 |           
     0       1
o3 : Complex
i4 : --  A simplicial complex realization of the Klein bottle 
as described in Armstrong's book "Basic Topology"
     L = abstractSimplicialComplex {{1,2,7}, {7,8,2}, {4,7,8}, 
     {4,8,5}, {1,4,5}, {1,5,2}, {2,8,3}, {8,3,9}, 
     {5,8,9}, {5,9,6}, {2,5,6},  {2,6,3}, {3,9,1}, {9,1,4}, 
     {6,9,4},{6,7,4}, {3,6,7}, {3,7,1}};
i5 : prune HH simplicialChainComplex L
       1
o5 = ZZ  <-- cokernel | 2 |
                      | 0 |
     0       1
o5 : Complex
\end{verbatim}
\end{scriptsize}
\end{example}

The design methodology utilized here, to implement a framework for abstract simplicial complexes within {\tt AbstractSimplicialComplexes.m2}, is based on graded lists.  Compared with the approach that is employed within the package {\tt SimplicialComplexes.m2}, there are a number of merits for the approach that we have adopted here.  

As a most basic example, the user interface for inputing the data of a simplicial complex within the package {\tt AbstractSimplicialComplexes.m2} is less cumbersome compared to that of {\tt SimplicialComplexes.m2}.  Indeed, when working within the package {\tt AbstractSimplicialComplexes.m2} there is no need to introduce a polynomial ring and then represent faces as products of the relevant variables.  Instead faces are defined simply as subsets of $[n] := \{1,\dots, n\} $ directly.  Further, the associated simplicial chain complexes (both reduced and non-reduced) are defined over $\ZZ$ by default.  By extension of scalars one obtains the reduced and non-reduced chain complexes with coefficients in other rings.  

As another example, our approach here facilitates a direct construction of random simplicial complexes.  This is an additional novel aspect to the present work which is of an independent interest.  It builds, for example, on the approach from \cite{Kahle:et:al:2020}.   Such methods to construct random simplicial complexes are currently not  implemented, directly, within the package {\tt SimplicialComplexes.m2}.  However, complementary constructions of random monomial ideals are implemented within the package {\tt RandomIdeals.m2}.  In either case, our approach here provides a new point of departure for such themes. 

In terms of computational benchmark comparable capabilities between {\tt AbstractSimplicialComplexes.m2} and {\tt SimplicialComplexes.m2},  their respective computable thresholds appear comparable (again depending on the particular applications that one has in mind).  Also, it is expected that the package {\tt AbstractSimplicialComplexes.m2} will see enhanced computational efficiencies once the main features of the package {\tt Complexes.m2}, \cite{Stillman:Smith:Complexes.m2}, on which the present package  {\tt AbstractSimplicialComplexes.m2} builds upon, become part of the core {\tt Macaulay2} infrastructure.    

To get a sense for the most basic computational comparable data between {\tt AbstractSimplicialComplexes.m2} and {\tt SimplicialComplexes.m2}, Example \ref{abstract:simplicial:complexes:SR:simplicialComplexes} below is helpful.  In considering this example, note that the methods used in {\tt SimplicialComplexes.m2} to construct the maps in the simplicial chain complexes that it produces ultimately rely on methods that are not internal to the package itself.   Rather they rely on methods that are internal to the core {\tt Macaulay2} infrastructure.  In doing so, some added efficiency, in terms of construction of simplicial chain complexes, is achieved.   

By contrast, in {\tt AbstractSimplicialComplexes.m2} the simplicial chain complexes are constructed directly within the package itself.  Aside from these computational considerations, in terms of construction of simplicial chain complexes, Example \ref{abstract:simplicial:complexes:SR:simplicialComplexes} suggests that one should expect there to be no major computational efficiency differences in terms of calculating the homology of such chain complexes and similar kinds of homological constructions.
 
\begin{example}[Comparing calculations on the $12$-simplex] \label{abstract:simplicial:complexes:SR:simplicialComplexes}
The following code snippet illustrates some of the most basic computational comparisons between  the respective packages {\tt AbstractSimplicialComplexes.m2} and {\tt SimplicialComplexes.m2}.  The conclusion is that, aside from apparent differences in creating the simplicial chain complexes themselves, compare {\tt i3} and {\tt i4} respectively, their respective comparable computational thresholds are similar (compare {\tt i4} and {\tt i5} respectively with {\tt i5} and {\tt i6}).  This is not surprising.

In the following code snippet, the the homology of the $12$-simplex is calculated using the methods of the package {\tt AbstractSimplicialComplexes.m2}.
\begin{scriptsize}
\begin{verbatim}
nathangrieve@Nathans-Air-3 ~ % M2
Macaulay2, version 1.24.11
with packages: ConwayPolynomials, Elimination, IntegralClosure, InverseSystems,
               Isomorphism, LLLBases, MinimalPrimes, OnlineLookup,
               PackageCitations, Polyhedra, PrimaryDecomposition, ReesAlgebra,
               Saturation, TangentCone, Truncations, Varieties
i1 : needsPackage"AbstractSimplicialComplexes";
i2 : time K = abstractSimplicialComplex(12);
 -- used 0.038916s (cpu); 0.0389109s (thread); 0s (gc)
i3 : time k = reducedSimplicialChainComplex(K);
 -- used 9.78308s (cpu); 7.7467s (thread); 0s (gc)
i4 : time h = HH k;
 -- used 38.1503s (cpu); 21.9464s (thread); 0s (gc)
i5 : time apply(13,i-> i => rank(h_i));
 -- used 102.599s (cpu); 50.6054s (thread); 0s (gc)
  \end{verbatim}
 \end{scriptsize}
The code snippet below reproduces these above calculations using the package {\tt SimplicialComplexes.m2}.
\begin{scriptsize}
\begin{verbatim}
nathangrieve@Nathans-Air-3 ~ % M2
Macaulay2, version 1.24.11
with packages: ConwayPolynomials, Elimination, IntegralClosure, InverseSystems,
               Isomorphism, LLLBases, MinimalPrimes, OnlineLookup,
               PackageCitations, Polyhedra, PrimaryDecomposition, ReesAlgebra,
               Saturation, TangentCone, Truncations, Varieties
i1 : needsPackage"SimplicialComplexes";
i2 : R = ZZ[x_1..x_12];
i3 : K = simplexComplex(11,R);
i4 : time k = complex K;
 -- used 0.104141s (cpu); 0.104134s (thread); 0s (gc)
i5 : time h = HH k;
 -- used 40.9266s (cpu); 22.6574s (thread); 0s (gc)
i6 : time apply(13,i-> i => rank(h_i));
 -- used 90.7713s (cpu); 47.5847s (thread); 0s (gc)
\end{verbatim}
\end{scriptsize}
\end{example}

Returning to general considerations, as some additional features, the package {\tt AbstractSimplicialComplexes.m2} is able to produce the induced chain complex maps that arise from nested pairs of simplicial complexes with one simplicial complex a subsimplicial complex of another.  The relevant methods for that purpose are  {\tt inducedSimplicialChainComplexMap} and {\tt reducedInducedSimplicialChainComplexMap}.   

As one possible direction for ongoing future development, which we do not pursue here, such methods may find applications towards creating an interface for the current package {\tt AbstractSimplicialComplexes.m2} to integrate with the package {\tt SpectralSequences.m2}.  Such an interface would also provide a simplification of the main constructors that are currently used within the package {\tt SpectralSequences.m2} to create filtered chain complexes arising from filtered simplicial complexes.  

Another direction to develop further is to expand on the existing infrastructure that we provide here so as to allow for applications that are in the direction of Topological Data Analysis (see for instance \cite{Carlsson:2009} and \cite{Carlsson:Vejdemo-Johansson}).  The design considerations that we have employed here will facilitate with this.  A first example in this direction is illustrated in Example \ref{random:point:cloud:eg}.  Further, the question of homological features for models of random simplicial complexes is quite attractive.  In Example \ref{random:simplicial:complex:models} we indicate one such model which arises in the work \cite{Kahle:et:al:2020} and the references therein.

The remainder of this article is organized in the following way.  In Section \ref{brief:overview:section} we give a brief overview of the {\tt AbstractSimplicialComplexes.m2} package and its syntax.  In Section \ref{random:simplicial:complexes} we discuss selected illustrative examples that pertain to homological calculations with random simplicial complexes.  In Section \ref{math:prelims} we provide a brief discussion of the relevant mathematical concepts that underly what we do here.  It also serves to fix the conventions and notation that are used implicitly within the package itself.  

\subsection*{Acknowledgements}
The author thanks the Natural Sciences and Engineering Research Council of Canada for their support through his grants DGECR-2021-00218 and RGPIN-2021-03821. 
It is the author's pleasure to thank colleagues for their interest and discussions on related topics.  The author also thanks Mike Stillman and  Greg Smith for helpful feedback in regards to the present work and also for corresponding to him about their package  \cite{Stillman:Smith:Complexes.m2}.  Finally, the author especially thanks an anonymous referee for their careful reading, testing of the code and thoughtful suggestions. 

\section{Selected Mathematical Preliminaries}\label{math:prelims}

\subsection{Abstract simplicial complexes}\label{simplicial:complex:prelims}

For our purposes here, by an \emph{abstract simplicial complex} $K$ on the vertex set 
$$[n] := \{1,\dots,n\}$$ 
we mean a collection of subsets 
$$\sigma \subseteq [n]$$ 
that is closed under taking subsets.  So in particular, if 
$$\tau \subseteq \sigma \subseteq [n] \text{ 
and }\sigma \in K \text{ then } \tau \in K \text{.}$$  

\subsection{Reduced and non-reduced simplicial chain complexes}\label{simplicial:chain:complex:prelims}

Fixing a commutative ring $R$, our conventions about reduced and non-reduced simplicial chain complexes with coefficients in $R$ are similar to those which can be found, for example, in \cite{CCA}.

Let $K$ be an abstract simplicial complex on the vertex set $[n]$.  Recall that 
$$
\dim K := \max_{\sigma \in K} \# \sigma -1 \text{.}
$$
Define $K$'s \emph{reduced simplicial chain complex} with coefficients in $R$ to be the chain complex
$$
0 \leftarrow \mathcal{C}_{-1}(K;R) \leftarrow \dots \leftarrow \mathcal{C}_{i-1}(K;R) \xleftarrow{\partial_i} \mathcal{C}_i(K;R) \leftarrow \dots \colon \mathcal{C}^{\mathrm{red}}_\bullet(K;R) \text{;}
$$
here
$$
\mathcal{C}_i(K;R) :=
\begin{cases}
\bigoplus\limits_{\substack{
\sigma \in K \\
\# \sigma = i + 1}} R \mathbf{e}_{\sigma} & \text{ for $-1 \leq i \leq \dim K$} \\
0 & \text{ otherwise} 
\end{cases}
$$
$$
\partial_i(\mathbf{e}_{\sigma}) := \sum_{j =0}^i (-1)^j \mathbf{e}_{\{\ell_0 < \dots < \hat{\ell}_j < \dots < \ell_i \}}  \text{ for $0 \leq i \leq \dim K$} 
$$
and  
$$\sigma = \{\ell_0 < \dots < \ell_i\} \in K\text{.}$$
(The notation $\{\ell_0 < \dots < \hat{\ell}_j < \dots < \ell_i \}$ means $\{\ell_0 < \dots < \ell_i\}$ with $\ell_j$ removed.)

The \emph{non-reduced simplicial chain complex} with coefficients in $R$ is defined to be the chain complex
$$
0 \leftarrow \mathcal{C}_{0}(K;R) \leftarrow \dots \leftarrow \mathcal{C}_{i-1}(K;R) \xleftarrow{\partial_i} \mathcal{C}_i(K;R) \leftarrow \dots \colon \mathcal{C}^{}_\bullet(K;R) \text{.}
$$
Respectively, denote the homology modules of these respective chain complexes by $H_i^{\mathrm{red}}(K;R)$ and $H_i(K;R)$ for $i \in \ZZ$.  So in particular
$$
H_i(K;R) = \ker \partial_i / \operatorname{image} \partial_{i+1} \text{ for $i \geq 0$}
$$
whereas 
$$
H_i^{\operatorname{red}}(K;R) = \ker \partial_i / \operatorname{image} \partial_{i+1} \text{ for $i \geq -1$.}
$$

\section{Brief overview of the package}\label{brief:overview:section}

In the package {\tt AbstractSimplicialComplexes.m2} abstract simplicial complexes are implemented by the type {\tt AbstractSimplicialComplexes}.  Users who wish to use the package {\tt AbstractSimplicialComplexes.m2} to study simplicial complexes on vertex sets different from $[n]$ should first fix a suitable order preserving bijection which is compatible with the standard lexicographic ordering.

Example \ref{simplicial:complex:constructors} below illustrates the main ways to construct simplicial complexes within the package {\tt AbstractSimplicialComplexes.m2}.

\begin{example}[Constructing abstract simplicial complexes]\label{simplicial:complex:constructors}
Within the package {\tt AbstractSimplicialComplexes.m2}, simplicial complexes can be constructed using the method {\tt abstractSimplicialComplexes}.  As input this method takes either:
\begin{itemize}
\item{an integer $n$ in which case the output is the $n$-simplex; or}
\item{a collection of subsets of $[n]$ in which case the output is the simplicial complex that is generated by these subsets.}
\end{itemize}

As another consideration to be aware of, if  $m \leq n$ and if $K$ is a simplicial complex with vertices supported on $[m]$ then automatically $K$ is considered as a subsimplicial complex of the $n$-simplex on $[n]$.  

Such features are  illustrated via the following code snippet.  
(More details about reduced and non-reduced simplicial chain complexes are given in Example \ref{eg:constructing:chain:complexes}.)

\begin{scriptsize}
\begin{verbatim}
nathangrieve@Nathans-Air-3 ~ % M2
Macaulay2, version 1.24.11
with packages: ConwayPolynomials, Elimination, IntegralClosure, InverseSystems,
               Isomorphism, LLLBases, MinimalPrimes, OnlineLookup,
               PackageCitations, Polyhedra, PrimaryDecomposition, ReesAlgebra,
               Saturation, TangentCone, Truncations, Varieties
i1 : needsPackage"AbstractSimplicialComplexes"
o1 = AbstractSimplicialComplexes
o1 : Package
i2 : --- make the simplicial complex on [5] 
     --- generated by {1,2,3,4},{1,3,4},{2,5},{2,4,5}
     K = abstractSimplicialComplex({{1,2,3,4},{1,3,4},{2,5},{2,4,5}});
i3 : -- vertices
     K_0
o3 = {{1}, {2}, {3}, {4}, {5}}
o3 : List
i4 : -- facets
     abstractSimplicialComplexFacets K
o4 = {{2, 4, 5}, {1, 2, 3, 4}}
o4 : List
i5 : -- make the simplex on [6]
     L = abstractSimplicialComplex(6);
i6 : -- we regard K as a subsimplicial complex of L
     -- the induced simplicial chain complex map is then
     f = inducedSimplicialChainComplexMap(L,K)

           6                      5
o6 = 0 : ZZ  <----------------- ZZ  : 0
                | 1 0 0 0 0 |
                | 0 1 0 0 0 |
                | 0 0 1 0 0 |
                | 0 0 0 1 0 |
                | 0 0 0 0 1 |
                | 0 0 0 0 0 |

           15                            8
     1 : ZZ   <----------------------- ZZ  : 1
                 | 1 0 0 0 0 0 0 0 |
                 | 0 1 0 0 0 0 0 0 |
                 | 0 0 1 0 0 0 0 0 |
                 | 0 0 0 0 0 0 0 0 |
                 | 0 0 0 0 0 0 0 0 |
                 | 0 0 0 1 0 0 0 0 |
                 | 0 0 0 0 1 0 0 0 |
                 | 0 0 0 0 0 1 0 0 |
                 | 0 0 0 0 0 0 0 0 |
                 | 0 0 0 0 0 0 1 0 |
                 | 0 0 0 0 0 0 0 0 |
                 | 0 0 0 0 0 0 0 0 |
                 | 0 0 0 0 0 0 0 1 |
                 | 0 0 0 0 0 0 0 0 |
                 | 0 0 0 0 0 0 0 0 |

           20                      5
     2 : ZZ   <----------------- ZZ  : 2
                 | 1 0 0 0 0 |
                 | 0 1 0 0 0 |
                 | 0 0 0 0 0 |
                 | 0 0 0 0 0 |
                 | 0 0 1 0 0 |
                 | 0 0 0 0 0 |
                 | 0 0 0 0 0 |
                 | 0 0 0 0 0 |
                 | 0 0 0 0 0 |
                 | 0 0 0 0 0 |
                 | 0 0 0 1 0 |
                 | 0 0 0 0 0 |
                 | 0 0 0 0 0 |
                 | 0 0 0 0 1 |
                 | 0 0 0 0 0 |
                 | 0 0 0 0 0 |
                 | 0 0 0 0 0 |
                 | 0 0 0 0 0 |
                 | 0 0 0 0 0 |
                 | 0 0 0 0 0 |

           15              1
     3 : ZZ   <--------- ZZ  : 3
                 | 1 |
                 | 0 |
                 | 0 |
                 | 0 |
                 | 0 |
                 | 0 |
                 | 0 |
                 | 0 |
                 | 0 |
                 | 0 |
                 | 0 |
                 | 0 |
                 | 0 |
                 | 0 |
                 | 0 |
o6 : ComplexMap
i7 : isWellDefined f
o7 = true
\end{verbatim}
\end{scriptsize}
\end{example}

\begin{example}[How to construct reduced and non-reduced simplicial chain complexes]\label{eg:constructing:chain:complexes}
There are two simplicial chain complex constructors.  The following code snippets illustrate the key points.

\begin{scriptsize}
\begin{verbatim}
i2 : K = abstractSimplicialComplex({{1,2,3,4}, {2,3,5}, {1,5}});
i3 : k = simplicialChainComplex K
       5       9       5       1
o3 = ZZ  <-- ZZ  <-- ZZ  <-- ZZ         
     0       1       2       3
o3 : Complex
i4 : k.dd
           5                                       9
o4 = 0 : ZZ  <---------------------------------- ZZ  : 1
                | -1 -1 -1 -1 0  0  0  0  0  |
                | 1  0  0  0  -1 -1 -1 0  0  |
                | 0  1  0  0  1  0  0  -1 -1 |
                | 0  0  1  0  0  1  0  1  0  |
                | 0  0  0  1  0  0  1  0  1  |
           9                           5
     1 : ZZ  <---------------------- ZZ  : 2
                | 1  1  0  0  0  |
                | -1 0  1  0  0  |
                | 0  -1 -1 0  0  |
                | 0  0  0  0  0  |
                | 1  0  0  1  1  |
                | 0  1  0  -1 0  |
                | 0  0  0  0  -1 |
                | 0  0  1  1  0  |
                | 0  0  0  0  1  |
           5               1
     2 : ZZ  <---------- ZZ  : 3
                | -1 |
                | 1  |
                | -1 |
                | 1  |
                | 0  |
o4 : ComplexMap
i5 : kRed = reducedSimplicialChainComplex K
       1       5       9       5       1
o5 = ZZ  <-- ZZ  <-- ZZ  <-- ZZ  <-- ZZ
     -1      0       1       2       3
o5 : Complex
i6 : kRed.dd
            1                      5
o6 = -1 : ZZ  <----------------- ZZ  : 0
                 | 1 1 1 1 1 |
           5                                       9
     0 : ZZ  <---------------------------------- ZZ  : 1
                | -1 -1 -1 -1 0  0  0  0  0  |
                | 1  0  0  0  -1 -1 -1 0  0  |
                | 0  1  0  0  1  0  0  -1 -1 |
                | 0  0  1  0  0  1  0  1  0  |
                | 0  0  0  1  0  0  1  0  1  |
           9                           5
     1 : ZZ  <---------------------- ZZ  : 2
                | 1  1  0  0  0  |
                | -1 0  1  0  0  |
                | 0  -1 -1 0  0  |
                | 0  0  0  0  0  |
                | 1  0  0  1  1  |
                | 0  1  0  -1 0  |
                | 0  0  0  0  -1 |
                | 0  0  1  1  0  |
                | 0  0  0  0  1  |
           5               1
     2 : ZZ  <---------- ZZ  : 3
                | -1 |
                | 1  |
                | -1 |
                | 1  |
                | 0  |
o6 : ComplexMap
\end{verbatim}
\end{scriptsize}

\end{example}

In Example \ref{sub:simplicial:complex}, we illustrate how to construct the induced chain complex maps that arise via subsimplicial complexes.

\begin{example}[Induced chain complex maps]\label{sub:simplicial:complex}
Given a subsimplicial complex $L \subseteq K$ there is support for computing the induced maps of the respective simplicial and reduced simplicial chain complexes. The following code snippets illustrate the key points.

\begin{scriptsize}
\begin{verbatim}
i2 : K = abstractSimplicialComplex({{1,2},{3}})
o2 = AbstractSimplicialComplex{-1 => {{}}          }
                               0 => {{1}, {2}, {3}}
                               1 => {{1, 2}}

o2 : AbstractSimplicialComplex
i3 : J = ambientAbstractSimplicialComplex(K)
o3 = AbstractSimplicialComplex{-1 => {{}}                   }
                               0 => {{1}, {2}, {3}}
                               1 => {{1, 2}, {1, 3}, {2, 3}}
                               2 => {{1, 2, 3}}
o3 : AbstractSimplicialComplex
i4 : f = inducedSimplicialChainComplexMap(J,K)
           3                  3
o4 = 0 : ZZ  <------------- ZZ  : 0
                | 1 0 0 |
                | 0 1 0 |
                | 0 0 1 |
           3              1
     1 : ZZ  <--------- ZZ  : 1
                | 1 |
                | 0 |
                | 0 |
o4 : ComplexMap
i5 : isWellDefined f
o5 = true
\end{verbatim}
\end{scriptsize}
\end{example}

\section{Selected discussion about random simplicial complexes}\label{random:simplicial:complexes}

\begin{example}[Models for random simplicial complexes]\label{random:simplicial:complex:models}
Within the package {\tt AbstractSimplicialComplexes.m2} we provide three ways of producing random simplicial complexes.  

The first method produces a random simplicial complex with support on $[n]$.  The second method produces a \emph{random simplicial complex} with support on $[n]$ and having $r$-skeleton.  The final method produces a random simplicial complex $Y_d(n,m)$ which has vertex set $[n]$, complete $(d-1)$ skeleton and exactly $m$ dimension $d$ faces chosen at random from all $\binom{\binom{n}{d+1}}{m}$ possibilities.

Such random simplicial complexes appear in lots of different contexts.  As some examples we refer to the article \cite{Kahle:et:al:2020} and the references therein. 

The following code snippets illustrate how to work with such methods using the package {\tt AbstractSimplicialComplexes.m2}.

\begin{scriptsize}
\begin{verbatim}
i2 : setRandomSeed(currentTime());
i3 : randomAbstractSimplicialComplex(5)
o3 = AbstractSimplicialComplex{-1 => {{}}     }
                               0 => {{1}, {3}}
                               1 => {{1, 3}}
o3 : AbstractSimplicialComplex
i4 : randomAbstractSimplicialComplex(5,3)
o4 = AbstractSimplicialComplex{-1 => {{}}                   }
                               0 => {{2}, {3}, {4}}
                               1 => {{2, 3}, {2, 4}, {3, 4}}
                               2 => {{2, 3, 4}}
o4 : AbstractSimplicialComplex
i5 : tally(for i from 1 to 10000 list (facets randomAbstractSimplicialComplex(5,3,2)))
o5 = Tally{{{1, 2, 3}, {1, 2, 4}, {2, 4, 5}} => 391 }
           {{1, 2, 3}, {1, 3, 4}, {3, 4, 5}} => 1338
           {{1, 2, 3}, {1, 3, 5}, {2, 3, 4}} => 1293
           {{1, 2, 3}, {1, 4, 5}, {2, 4, 5}} => 403
           {{1, 2, 4}, {1, 2, 5}} => 1710
           {{1, 4, 5}, {2, 4, 5}} => 1774
           {{2, 3, 4}, {2, 3, 5}, {3, 4, 5}} => 1557
           {{2, 3, 4}, {2, 4, 5}, {3, 4, 5}} => 1534
o5 : Tally
\end{verbatim}
\end{scriptsize}

\end{example}

\begin{example}[Random point clouds and Vietoris-Rips complexes]\label{random:point:cloud:eg}
In what follows we illustrate a collection of homological calculations that can be performed on random Vietoris-Rips complexes.  

Indeed, motivated by the perspective of Topological Data Analysis, see for example \cite{Carlsson:Vejdemo-Johansson}, here we model a dimension $n$ \emph{random point cloud distance matrix} as a random  $n \times n$ upper triangular matrix with random real entries above the main diagonal.  

In doing so the $(i,j)$-entry represents the distance $d(x_i,x_j)$ between the points $x_i$ and $x_j$; the corresponding \emph{Vietoris-Rips complex} with vertex set supported on $[n]$ and depending on parameter $\epsilon > 0$ is the simplicial complex that has $k$-faces those cardinality $k+1$ subsets
$$\{i_0 < \dots < i_k\} \subseteq [n]$$ 
which are such that 
$$d(x_{i_j},x_{i_{\ell}}) \leq \epsilon $$ 
$$\text{ for all $i_j, i_{\ell} \in \{i_0 < \dots < i_k\}$.}$$
We refer to of the text \cite[p.~104]{Carlsson:Vejdemo-Johansson} and the survey article \cite[Sec.~2.2]{Carlsson:2009} for more details.
\begin{scriptsize}
\begin{verbatim}

i2 : -- A random such VR complex on [n] and depending on a given
     -- parameter e can be modelled in the following way.
     randomVRcomplex := (n,e) -> (
     -- return a random VR-complex supported on [n] and depending on parameter e --
         L := for i from 1 to n list i;
         setRandomSeed(currentTime());
         M := fillMatrix(mutableMatrix(RR,n,n), UpperTriangular => true);
         myFaces := select(subsets(L), i-> all(apply(subsets(i,2), 
         j-> M_(j#0-1,j#1-1) <= e)));
         K = abstractSimplicialComplex myFaces;
         return K);
i3 : -- So for example the facets of a random VR complex on the vertex set [10] 
     -- and depending on a parameter e = .4 is described as
     setRandomSeed(currentTime());
i4 : K = randomVRcomplex(10,.4);
i5 : facets K
o5 = {{5, 8}, {2, 4, 9}, {3, 4, 6}, {4, 6, 9}, {4, 7, 10}, {6, 8, 9}, 
{7, 8, 9}, {1, 4, 7, 9}, 
{2, 3, 4, 10}, {3, 4, 5, 10}}
o5 : List
i6 : prune HH reducedSimplicialChainComplex K
       1
o6 = ZZ
     1
o6 : Complex
\end{verbatim}
\end{scriptsize}
\end{example}

\providecommand{\bysame}{\leavevmode\hbox to3em{\hrulefill}\thinspace}
\providecommand{\MR}{\relax\ifhmode\unskip\space\fi MR }
\providecommand{\MRhref}[2]{%
  \href{http://www.ams.org/mathscinet-getitem?mr=#1}{#2}
}
\providecommand{\href}[2]{#2}


\end{document}